\documentclass [12pt]{article}
\usepackage {amssymb}
\usepackage {amsmath}

\begin{document}

\begin{center}
\textbf{A~CHARACTERIZATION OF~ROOT CLASSES OF~GROUPS}
\end{center}

\begin{center}
\textsc{E.~V.~Sokolov}
\end{center}

\begin{abstract}
We prove that a~class of~groups is root in~a~sense of~K.~W.~Gruenberg if, and~only if, it is closed under subgroups and~Cartesian wreath products. Using this result we obtain a~condition which is sufficient for the~generalized free product of~two nilpotent groups to~be residual solvable.
\end{abstract}

Recall that a~class of~groups $\mathcal{C}$ is called root if~it satisfies the~following conditions:

1)\hspace{0.5em}an~arbitrary subgroup of~a \hbox{$\mathcal{C}$-}group is also a~\hbox{$\mathcal{C}$-}group;

2)\hspace{0.5em}the~direct product of~any two \hbox{$\mathcal{C}$-}groups belongs to~$\mathcal{C}$;

3)\hspace{0.5em}the~Gruenberg condition: for any group~$X$ and~for any subnormal sequence $Z \leqslant Y \leqslant X$ with~factors in~$\mathcal{C}$, there exists a~normal subgroup~$T$ of~$X$ such that $T \leqslant Z$ and~$X/T \in \mathcal{C}$.

It is easy to~see that the~classes of~all finite groups, of~all finite $p$-groups, and~of~all solvable groups are root. The class of~all nilpotent groups fails to~be root because it does not satisfy the~Gruenberg condition.

The notion of~the root class was introduced by~K.~W.~Gruenberg~\cite{li01} in~connection with~studying of~residual properties of~solvable groups. D.~N.~Azarov and~D.~Tieudjo~\cite{li02} proved that every free group is residually a~\hbox{$\mathcal{C}$-}group for any nontrivial (i.~e. containing at least one non-unit group) root class of~groups. It follows from this statement and~the~results of~K.~W.~Gruenberg~\cite{li01} that, for every nontrivial root class of~groups $\mathcal{C}$, the~free product of~an arbitrary number of~residually \hbox{$\mathcal{C}$-}groups is itself residually a~\hbox{$\mathcal{C}$-}group. The property `to~be residually a~\hbox{$\mathcal{C}$-}group', where $\mathcal{C}$~is an~arbitrary root class of~groups, was studied in~\cite{li03} and~\cite{li04} in~relation to~generalized free products and~HNN-extensions. The results of~these papers show that many conditions known to~be necessary or~sufficient for a~given group to~be residually a~\hbox{$\mathcal{C}$-}group for some concrete root class~$\mathcal{C}$ remain true as well in~the~case when $\mathcal{C}$~is an~arbitrary root class of~groups.

The~definition given by~K.~W.~Gruenberg doesn't allow to~describe all root classes of~groups clearly. The~main aim of~this paper is to~give another, more simple, characterization of~root classes.

It is easy to~see that the~second condition in~the~definition of~the root class follows from the~first and~the~third and, thus, is excessive. As~for the~Gruenberg condition, it can be replaced by~a~more clear criterion as~the~next theorem shows.

\medskip

\textbf{Theorem~1.}~Let $\mathcal{C}$ be a~class of~groups closed under taking subgroups. Then the~following statements are equivalent:

1.\hspace{0.5em}$\mathcal{C}$ satisfies the~Gruenberg condition (and, hence, is root).

2.\hspace{0.5em}$\mathcal{C}$ is closed under Cartesian wreath products.

3.\hspace{0.5em}$\mathcal{C}$ is closed under extensions and, for any two groups $X, Y \in \mathcal{C}$, contains the~Cartesian product $\prod_{y \in Y}X_{y}$, where $X_{y}$~is an~isomorphic copy of~$X$ for every $y \in Y$.

\medskip

The next three statements follow immediately from Theorem~1.

\medskip

\textbf{Corollary~1} \cite{li05}. If a~class of~groups consists of~the only finite groups, then it is root if, and~only if, it is closed under subgroups and~extensions.~$\square$

\medskip

\textbf{Corollary~2}. The~intersection of~any two root classes of~groups is again a~root class of~groups.~$\square$

\medskip

\textbf{Corollary~3.} If~$\mathcal{C}$ is a~root class of~groups, then the~class $\mathcal{C}_{0}$ of~all torsion-free \hbox{$\mathcal{C}$-}groups is root too.~$\square$

\medskip

The theorem given below serves as an~example of~application of~Theorem~1 to~studying of~residual properties of~generalized free products of~groups.

\medskip

\textbf{Theorem~2.} Let $\mathcal{C}$ be a~nontrivial root class of~groups closed under taking quotient groups. Let also $G = \langle A * B;\ H = K,\ \varphi \rangle$ be the~free product of~nilpotent \hbox{$\mathcal{C}$-}groups~$A$ and~$B$ with~subgroups $H \leqslant A$ and~$K \leqslant B$ amalgamated according to~an~isomorphism~$\varphi$. Suppose that $A$ and~$B$ possess such central series
$$
1 = A_{0} \leqslant A_{1} \leqslant \ldots \leqslant A_{n} = A,\ 
1 = B_{0} \leqslant B_{1} \leqslant \ldots \leqslant B_{n} = B,
$$
that $(A_{i} \cap H)\varphi = B_{i} \cap K$ for every $i \in \{1, 2, \ldots, n\}$. Then there exists a~homomorphism of~$G$ onto~a~solvable \hbox{$\mathcal{C}$-}group which is injective on~$A$ and~$B$. In~particular, $G$ is residually a~\hbox{$\mathcal{CS}$-}group, where $\mathcal{CS}$ is the~class of~all solvable \hbox{$\mathcal{C}$-}groups.

\medskip

We note that Theorem~2 strengthens and~generalizes Theorem~8 from \cite{li06}, which states that $G$ is poly-(residually solvable).

\medskip

\textit{Proof of~Theorem~1.}

$1 \Rightarrow 2$. Let $X$, $Y$ be arbitrary \hbox{$\mathcal{C}$-}groups, and~let $B$ be the~Cartesian product of~isomorphic copies of~$X$ indexed by~the~elements of~$Y$ (i.~e. $B$ is the~set of~all functions mapping $Y$ into~$X$ with~the~coordinate-wise multiplication). Let also $W = X \wr Y$ be the~Cartesian wreath product of~$X$ and~$Y$. We need to~show that $W \in \mathcal{C}$.

Recall that $W$ is the~set $Y \cdot B$ with~the~operation defined by~the~rule $y_{1}b_{1}y_{2}b_{2} = y_{1}y_{2}b^{y_2}b_{2}$, where $b^{y_2} \in B$ is the~function mapping $y$ to~$y_{2}y$ for every $y \in Y$. From this definition it follows that $B$ is normal in~$W$ and~$W/B \cong Y$.

Let $A = \{b \in B \mid b(1) = 1\}$. Then $A$ is normal in~$B$ and~$B/A \cong X$. Since $W/B \cong Y$, then $A \leqslant B \leqslant W$ is a~subnormal sequence with~\hbox{$\mathcal{C}$-}factors, and, by~the~Gruenberg condition, there exists a~normal subgroup~$T$ of~$W$ such that $T \leqslant A$ and~$W/T \in \mathcal{C}$.

As~$T$ is normal in~$W$, it is contained in~the~subgroup
$$
A^{y} = \{b \in B \mid b(y) = 1\}
$$
for any $y \in Y$. But~$\bigcap_{y \in Y}A^{y} = 1$, hence, $T = 1$ and~$W \in \mathcal{C}$ as~required.

$2 \Rightarrow 3$. Let $X$, $Y$ be arbitrary \hbox{$\mathcal{C}$-}groups, $W = X \wr Y$, and~$B = \prod_{y \in Y}X_{y}$, where $X_{y}$ is an~isomorphic copy of~$X$ for every $y \in Y$. Then $W \in \mathcal{C}$, $B \leqslant W$, and~$B \in \mathcal{C}$ since $\mathcal{C}$ is closed under subgroups. Further, if~$Z$ is an~extension of~$X$ by~$Y$, then, by~the~theorem of~Frobenius, $Z$ embeds in~$W$ and,~therefore, belongs to~$\mathcal{C}$.

$3 \Rightarrow 1$. Let~$X$ be a~group, and~let $Z \leqslant Y \leqslant X$ be a~subnormal sequence with~\hbox{$\mathcal{C}$-}factors. We put $T = \bigcap_{s \in S}Z^{s}$, where $S$ is some system of~all cosets representatives of~$Y$ in~$X$, and~show that~$T$ is required.

It is obvious that~$T$ is a~normal subgroup of~$X$ lying in~$Z$. The quotient group $Y/T$ embeds in~the~Cartesian product $P$ of~the quotient groups $Y/Z^{s}$, $s \in S$, by~the~theorem of~Remak. Each of~groups $Y/Z^{s}$ is isomorphic to~the~\hbox{$\mathcal{C}$-}group~$Y/Z$. Therefore, $P \in \mathcal{C}$, and~$Y/T \in \mathcal{C}$ since $\mathcal{C}$ is closed under subgroups. Thus, $Y/T \in \mathcal{C}$, $X/Y \in \mathcal{C}$, and~$X/T \in \mathcal{C}$ because $\mathcal{C}$ is closed under extensions.~$\square$

\medskip

\textit{Proof of~Theorem~2} will use an~induction on~$n$.

If~$n = 1$, then there exists a~homomorphism of~$G$ onto~the~generalized direct product $P = \langle A \times B;\ H = K,\ \varphi \rangle$ continuing the~natural inclusions of~$A$ and~$B$, and~this homomorphism is required. Indeed, $P$ is isomorphic to~the~quotient group of~the direct product $A \times B$ by~the~subgroup $gp\{h(h\varphi)^{-1} \mid h \in H\}$. Since $\mathcal{C}$ is a~root class, then $A \times B \in \mathcal{C}$. It remains to~note that $\mathcal{C}$ is closed under quotient groups and~so $P \in \mathcal{C}$.

Let now $n > 1$, and~let $\bar\varphi\colon HA_{1}/A_{1} \to KB_{1}/B_{1}$ be a~map such that $(hA_{1})\bar\varphi = (h\varphi)B_{1}$ for every $h \in H$. It follows from the~equality $(A_{1} \cap H)\varphi = B_{1} \cap K$ that $\bar\varphi$ is a~correctly defined isomorphism of~subgroups. Therefore, we can consider the~generalized free product $\bar G = \langle \bar A * \bar B;\ \bar H = \bar K,\ \bar\varphi \rangle$, where $\bar A = A/A_{1}$, $\bar B = B/B_{1}$, $\bar H = HA_{1}/A_{1}$, and~$\bar K = KB_{1}/B_{1}$.

Since $\mathcal{C}$ is closed under quotient groups, $\bar A, \bar B \in \mathcal{C}$. It is easy to~see also that the~series
\begin{gather*}
1 = A_{1}/A_{1} \leqslant A_{2}/A_{1} \leqslant \ldots \leqslant A_{n}/A_{1} = \bar A,\\ 
1 = B_{1}/B_{1} \leqslant B_{2}/B_{1} \leqslant \ldots \leqslant B_{n}/B_{1} = \bar B
\end{gather*}
are $\bar\varphi$-compatible. Hence, $\bar G$ satisfies the~conditions of~the theorem and, by~induction hypothesis, there exists a~homomorphism of~this group onto~a~solvable \hbox{$\mathcal{C}$-}group~$Y$ which is injective on~$\bar A$ and~$\bar B$.

Since $\mathcal{C}$ is closed under subgroups, $A_{1}, B_{1} \in \mathcal{C}$. Therefore, the~group \hbox{$G_{1} = \langle A_{1} * B_{1};\ H_{1} = K_{1},\ \varphi_{1} \rangle$,} where $H_{1} = H \cap A_{1}$, $K_{1} = K \cap B_{1}$, and \hbox{$\varphi_{1} = \varphi\vert_{H_1}$,} satisfies the~conditions of~the theorem too. Again by~induction hypothesis, there exists a~homomorphism of~this group onto~a~solvable \hbox{$\mathcal{C}$-}group~$X$ which is injective on~$A_{1}$ and~$B_{1}$.

Now we use Lemma~2 from~\cite{li07}. By this lemma, there exists a~homomorphism~$\rho$ of~$G$ onto~$X \wr Y$ which is injective on~$A$ and~$B$. $X \wr Y$ is a~solvable group and, by~Theorem~1, belongs to~$\mathcal{C}$. Hence, $\rho$ is the~required homomorphism.

By the~condition of~the theorem, $\mathcal{C}$ contains at least one nontrivial group. All cyclic subgroups of~this group and~their quotient groups belong to~$\mathcal{C}$ too. Hence, $\mathcal{C}$ includes at least one cyclic group of~prime order, say~$p$. It is well known that every finite $p$-group possesses a~normal series with~the~factors of~order~$p$. Taking into~account that $\mathcal{C}$ is closed under extensions we conclude that all finite $p$-groups are contained in~$\mathcal{C}$.

Let $N = \ker\rho$. Since $N \cap A = N \cap B = 1$, then, by~the~theorem of~H.~Neumann~\cite{li08}, $N$ is free. As it is known, free groups are residually $p$-finite for any prime~$p$. Hence, $N$ is residually a~\hbox{$\mathcal{CS}$-}group.

Thus, $G$ is an~extension of~the residually \hbox{$\mathcal{CS}$-}group~$N$ by~the~\hbox{$\mathcal{CS}$-}group $G/N$. By Corollary~2, the~class of~all solvable \hbox{$\mathcal{C}$-}groups is root. Therefore, $G$ is a~residually \hbox{$\mathcal{CS}$-}group by~Lemma~1.5 from~\cite{li01}.~$\square$

\end{document}